\documentclass[12pt]{article}
\usepackage{amsfonts}
\usepackage{amsmath,amssymb,amsfonts}

\newtheorem{theorem}{Theorem}

\newtheorem{corollary}[theorem]{Corollary}
\newtheorem{definition}[theorem]{Definition}

\newtheorem{remarks}[theorem]{Remark}

\begin{document}

\title{A dynamical presentation of the better than nice metric on the disc}
\author{Daniela Emmanuele and Marcos Salvai\thanks{%
Partially supported by Conicet, Foncyt, Secyt Univ.\thinspace Nac.\thinspace
C\'{o}rdoba}}
\date{}
\maketitle

\begin{abstract}
This is a postprint of our paper \textsl{Force free Moebius motions of the
circle} (J. Geom. Symmetry Phys. 27 (2012) 59--65), which we hadn't uploaded
to arXiv previously. We would like to draw attention to the relationship
with the article \textsl{A geometry where everything is better than nice},
by Larry Bates and Peter Gibson (to appear in Proc. Amer. Math. Soc.).

In our note we treat thoroughly a simple particular case of two previous,
more substantial articles. We describe the force free Moebius motions of the
circle, that is, the geodesics of the Lie group $\mathcal{M}\simeq \mathrm{%
PSL}_{2}\left( \mathbb{R}\right) \simeq \mathrm{O}_{o}\left( 1,2\right) $ of
Moebius transformations of the circle, equipped with the Riemannian metric
given by the kinetic energy induced by the action. It turns up that $%
\mathcal{M}$ decomposes as a Riemannian product $S^{1}\times \Delta $, where 
$\Delta $ is the unit disc endowed with a certain metric, which we now
recognize as being essentially the one that is better than nice. Concerning
geodesics, we only give the differential equation for their trajectories
(using Clairaut's Theorem); we took pleasure in learning from the article by
Bates and Gibson that they are actually hypocycloids.
\end{abstract}

\noindent MSC 2010: 70\thinspace K\thinspace 25, 22\thinspace E\thinspace
70, 53\thinspace C\thinspace 22.

\smallskip


\section{Introduction}

In the spirit of the classical description of the force free motions of a
rigid body in Euclidean space using an invariant metric on $\mathrm{SO}%
\left( 3\right) $ \cite[Appendix~2]{arnold}, the second author defined in 
\cite{Salvai1} an appropriate metric on the Lorentz group $\mathrm{SO}%
_{o}\left( n+1,1\right) $ to study force free conformal motions of the
sphere ${\mathbb{S}}^{n}$, obtaining a few explicit ones (only through the
identity and those which can be described using the Lie structure of the
configuration space). In this note, in the particular case $n=1$, that is,
Moebius motions of the circle, we obtain all force free motions.

This is an example of a situation in which using concepts of Physics one can
state and solve a problem in Differential Geometry.

Notice that the canonical action of $\mathrm{PSL}\left( 2,\mathbb{R}\right) $
on $\mathbb{R}P^{1}\cong {\mathbb{S}}^{1}$ is equivalent to the action of
the group of Moebius transformations on the circle. Then, the results
presented here, up to a double covering, also extend the case $n=1$ of \cite%
{Salvai2}, where force free projective motions of the sphere ${\mathbb{S}}%
^{n}$ were studied.

This note, as well as \cite{Salvai1, Salvai2}, is weakly related with mass
transportation \cite{villani}. In our situation, the set of admitted mass
distributions is finite dimensional, and also the allowed transport maps are
very particular.

We would like to mention that the closure of $\mathcal{M}$ in the space of
all measurable functions of the circle into itself provides a
compactification of $\mathcal{M}$ as an open dense subset of the
three-sphere, with a dynamical meaning related to generalized flows (see 
\cite{Compactif}).

\subsection{Moebius motions of the circle}

Let ${\mathbb{S}}^{1}$ be the unit circle centered at zero in $\mathbb{C}$
with the usual metric and let $\mathcal{M}$ be the Lie group of Moebius
transformations of the circle, that is, the group of Moebius transformations
of the extended plane preserving the circle. It consists of maps of the form 
$cT_{\alpha },$ where $c\in {\mathbb{S}}^{1}$ and

\begin{equation}
T_{\alpha }\left( z\right) =\frac{z+\alpha }{1+\bar{\alpha}z}  \label{Talpha}
\end{equation}%
for $\alpha \in \mathbb{C}$, $\left\vert \alpha \right\vert <1$ and all $%
z\in {\mathbb{S}}^{1}.$ Although we are interested in the action of $%
\mathcal{M}$ on the circle, we recall that if the unit disc $\Delta =\left\{
z\in \mathbb{C}\nobreak\mskip2mu\mathpunct{}\nonscript\mkern-\thinmuskip{;} %
\mskip6muplus1mu\relax \left\vert z\right\vert <1\right\} $ carries the
canonical Poincar\'{e} metric of constant negative curvature $-1$ and $%
\alpha \neq 0,$ then $T_{\alpha }$ is the transvection translating the
geodesic with end points $\pm \alpha /\left\vert \alpha \right\vert $,
sending $0$ to $\alpha .$

A \emph{Moebius motion} of the circle is by definition a smooth curve in $%
\mathcal{M}$, thought of as a curve of diffeomorphisms of the circle.
(Throughout the paper, smooth means of class $C^{\infty }$.)

\bigskip%

In the next two subsections we recall, specialized for the circle, some
definitions and statements given in \cite{Salvai1} for conformal motions on
the $n$-dimensional sphere.

\subsection{The energy of Moebius motions of the circle}

Suppose that the circle has initially a homogeneous distribution of mass of
constant density $1$ and that the particles are allowed to move only in such
a way that two configurations differ in an element of $\mathcal{M}$. The
configuration space may be naturally identified with $\mathcal{M}$.

Let $\gamma :\left[ t_{0},t_{1}\right] \rightarrow \mathcal{M}$ be a Moebius
motion of ${\mathbb{S}}^{1}$. The total kinetic energy $E_{\gamma }\left(
t\right) $ of the motion $\gamma $ at the instant $t$ is given by 
\begin{equation}
E_{\gamma }\left( t\right) =\tfrac{1}{2}\int_{{\mathbb{S}}^{1}}\left\vert
v_{t}\left( q\right) \right\vert ^{2}\rho _{t}\left( q\right) \,{d}m\left(
q\right) \text{,}  \label{energy1}
\end{equation}%
where integration is taken with respect to the canonical volume form of ${%
\mathbb{S}}^{1}$ and, if $q=\gamma \left( t\right) \left( p\right) $ for $%
p\in {\mathbb{S}}^{1}$, then 
\begin{equation*}
v_{t}\left( q\right) =\left. \frac{{d}}{{d}s}\right\vert _{t}\gamma \left(
s\right) \left( p\right) \in T_{q}{\mathbb{S}}^{1}\text{,\ \ \ \ \ \ \ \ }%
\rho _{t}\left( q\right) =1/\det \left( {d}\gamma (t)_{p}\right)
\end{equation*}%
are the velocity of the particle $q$ and the density at $q$ at the instant $%
t $, respectively. Applying to (\ref{energy1}) the formula for change of
variables, one obtains 
\begin{equation}
E_\gamma \left( t\right) =\tfrac{1}{2}\int_{{\mathbb{S}}^{1}}\left\vert
\left. \frac{{d}}{{d}s}\right\vert _{t}\gamma \left( s\right) \left(
p\right) \right\vert ^{2}\,{d}m\left( p\right) \text{.}  \label{energy}
\end{equation}

The kinetic energy of $\gamma $ is defined by%
\begin{equation*}
E\left( \gamma \right) =\int_{t_{0}}^{t_{1}}E_{\gamma }\left( t\right) ~{d}t%
\text{. }
\end{equation*}%
The following definition is based on the principle of least action.

\medskip%

\noindent%

\begin{definition}
A smooth curve $\gamma $ in $\mathcal{M}$, thought of as a Moebius motion of 
${\mathbb{S}}^{1}$, is said to be \emph{force free }if it is a critical
point of the kinetic energy functional, that is, 
\begin{equation*}
\left. \frac{{d}}{{d}s}\right\vert _{0}E\left( \gamma _{s}\right) =0
\end{equation*}%
for any proper smooth variation $\gamma _{s}$ of $\gamma $ \emph{(}here $%
\gamma _{s}\left( t\right) =\Gamma \left( s,t\right) $, where $\Gamma
:\left( -\varepsilon ,\varepsilon \right) \times \left[ t_{0},t_{1}\right]
\rightarrow \mathcal{M}$ is a smooth map, with $\varepsilon >0$, $\Gamma
\left( 0,t\right) =\gamma \left( t\right) $ and $\Gamma \left(
s,t_{i}\right) =\gamma \left( t_{i}\right) $ for all $s\in \left(
-\varepsilon ,\varepsilon \right) $, $i=0,1$\emph{)}.
\end{definition}

\subsection{A Riemannian metric on the configuration space}

Given $g\in \mathcal{M}$ and $X\in T_{g}\mathcal{M}$, let us define the map $%
\widetilde{X}:{\mathbb{S}}^{1}\rightarrow T{\mathbb{S}}^{1}$ by 
\begin{equation}
\widetilde{X}(q)=\left. \frac{{d}}{{d}t}\right\vert _{0}\gamma \left(
t\right) \left( q\right) \in T_{g\left( q\right) }{\mathbb{S}}^{1}\text{,}
\label{Xtilde}
\end{equation}
where $\gamma $ is any smooth curve in $\mathcal{M}$ with $\gamma \left(
0\right) =g$ and $\dot{\gamma}\left( 0\right) =X$. The map $\widetilde{X}$
is well-defined and smooth and it is a vector field on ${\mathbb{S}}^{1}$ if
and only if $X\in T_{e}\mathcal{M}$. Moreover, 
\begin{equation}
X\mapsto \left\Vert X\right\Vert ^{2}=\frac{1}{2\pi }\int_{{\mathbb{S}}%
^{1}}\left\vert \widetilde{X}(q)\right\vert ^{2}\,{d}m\left( q\right)
\label{metricG}
\end{equation}
is a quadratic form on $T_{g}\mathcal{M}$ and gives a Riemannian metric on $%
\mathcal{M}$.

\begin{remarks}
\emph{a)} The fundamental property of the metric \emph{(\ref{metricG})} on $%
\mathcal{M}$ is that a curve $\gamma $ in $\mathcal{M}$ is a geodesic if and
only if \emph{(}thought of as a Moebius motion\emph{)} it is force free,
since by \emph{(\ref{metricG})} and \emph{(\ref{energy})}, $E_\gamma\left(
t\right) =\pi \left\Vert \dot{\gamma}\left( t\right) \right\Vert ^{2}$.

\medskip

\emph{b)} The metric on $\mathcal{M}$ is neither left nor right invariant,
since we saw in \emph{\cite{Salvai1}} that it is not even complete. \bigskip
\end{remarks}

\section{Force free Moebius motions of the circle}

The next theorem describes completely the geometry of $\mathcal{M}$ endowed
with the metric (\ref{metricG}) given by the kinetic energy. Recall from (%
\ref{Talpha}) that $T_{a}$ denotes the transvection associated with $\alpha $
and that $\Delta $ is the unit disc centered at zero in $\mathbb{C}$.

\begin{theorem}
\label{Th} Let ${d}s^{2}$ be the metric on the disc $\Delta $ given in polar
coordinates $\left( r,\theta \right) $ by 
\begin{equation}
{d}s^{2}=\frac{2\left( {d}r^{2}+r^{2}{d}\theta ^{2}\right) }{1-r^{2}}
\label{metricaDelta}
\end{equation}%
and consider on ${\mathbb{S}}^{1}\times \Delta $ the Riemannian product
metric, where ${\mathbb{S}}^{1}$ has length $2\pi $. Then the map 
\begin{equation*}
F:{\mathbb{S}}^{1}\times \Delta \rightarrow \mathcal{M}\text{,\ \ \ \ \ }%
F\left( u,\alpha \right) =uT_{\alpha }
\end{equation*}%
is an isometry.
\end{theorem}

\begin{remarks}
\emph{a)} Note that the metric \emph{(}\ref{metricaDelta}\emph{)}\ on $%
\Delta $ is not the canonical metric of constant negative curvature on $%
\Delta $. Indeed, the curvature function can be easily computed to be $%
K(r,\theta )=-1/\left( 1-r^{2}\right) $, in particular, it tends to $-\infty 
$ as $r\rightarrow 1^{-}$. Also, the metric on $\Delta $ is not complete,
since the inextendible ray $\left( 0,1\right) \backepsilon r\mapsto T_{r}$
has length $\pi /\sqrt{2}$, since $\left\Vert \frac{\partial }{\partial r}%
\right\Vert ^{2}=\frac{2}{1-r^{2}}$.

\smallskip

\emph{b)} In the higher dimensional situation \emph{\cite{Salvai1}} it is
proven that the group $\mathrm{SO}(n)$ \emph{(}with the metric induced from
the one given by the kinetic energy\emph{)} is totally geodesic in the group
of directly conformal transformations of ${\mathbb{S}}^{n}$, but the author
did not know whether this subgroup is a Riemannian factor, as it turned to
be for $n=1$. In the projective case \emph{\cite{Salvai2}}, $\mathrm{SO}(n)$
is not even totally geodesic.
\end{remarks}

\noindent\textbf{Proof of Theorem \ref{Th}.} Let ${\mathbb{S}}^{1}\subset 
\mathcal{M}$ be the subgroup of isometries of the circle. The torus ${%
\mathbb{S}}^{1}\times {\mathbb{S}}^{1}$ acts on $\mathcal{M}$ on the left by 
$\left( u,v\right) \cdot g=ug\bar{v}$, where $\left( ug\bar{v}\right) \left(
z\right) =ug\left( z\bar{v}\right) $ for any $z\in {\mathbb{S}}^{1}$. We
know from the higher dimensional cases in \cite{Salvai1} that this action is
by isometries of $\mathcal{M}$, provided that this group is endowed with the
metric (\ref{metricG}).

We fix $0<r<1$. By the torus symmetry just described, it suffices to verify
that ${d}F_{\left( 1,r\right) }:T_{\left( 1,r\right) }\left( {\mathbb{S}}%
^{1}\times \Delta \right) \rightarrow T_{F\left( 1,r\right) }\mathcal{M}$ is
a linear isometry. We put coordinates $t\mapsto {e}^{{i}t}$ on ${\mathbb{S}}%
^{1}$ and $\left( \rho ,\theta \right) \mapsto \rho {e}^{{i}\theta }$ on $%
\Delta $. We denote $\partial _{x}=\frac{{d}}{{d}x}$. Let $X,Y,Z$ be the
images under ${d}F_{\left( 1,r\right) }$ of $\partial _{t},\partial _{\rho
},\partial _{\theta }$, respectively. It suffices to show that $\left\{
X,Y,Z\right\} $ is an orthogonal basis of $T_{F\left( 1,r\right) }\mathcal{M}
$ with 
\begin{equation*}
\left\| X\right\| ^{2}=1\text{,\ \ \ \ \ \ \ }\Vert Y\Vert ^{2}=\frac{2}{%
1-r^{2}}\text{, \ \ \ \ \ \ }\Vert Z\Vert ^{2}=\frac{2r^{2}}{1-r^{2}}\text{.}
\end{equation*}

First, we compute $\widetilde{X}$, $\widetilde{Y}$ and $\widetilde{Z}$ by
their definition (\ref{Xtilde}). In each case, we take the curve $\gamma $
as the image under $F$ of the coordinate curves in ${\mathbb{S}}^{1}\times
\Delta $ through the point $\left( 1,r\right) $. We have 
\begin{eqnarray*}
\widetilde{X}(z) &=&\left. \frac{{d}}{{d}t}\right\vert _{0}F({e}^{{i}%
t},r)\left( z\right) =\left. \frac{{d}}{{d}t}\right\vert _{0}{e}^{{i}%
t}T_{r}(z)=\left. {i}{e}^{{i}t}T_{r}(z)\right\vert _{t=0}={i}T_{r}(z)={i}%
\frac{z+r}{1+rz} \\
&& \\
\widetilde{Y}\left( z\right) &=&\left. \frac{{d}}{{d}\rho }\right\vert
_{_{r}}F(1,\rho )\left( z\right) =\left. \frac{{d}}{{d}\rho }\right\vert
_{_{r}}T_{\rho }(z)=\left. \frac{{d}}{{d}\rho }\right\vert _{_{r}}\frac{%
z+\rho }{1+\rho z}=\frac{1-z^{2}}{(1+rz)^{2}} \\
&& \\
\widetilde{Z}\left( z\right) &=&\left. \frac{{d}}{{d}\theta }\right\vert
_{0}F(1,r{e}^{{i}\theta })\left( z\right) =\left. \frac{{d}}{{d}\theta }%
\right\vert _{0}T_{r{e}^{{i}\theta }}(z)=\left. \frac{{d}}{{d}\theta }%
\right\vert _{0}\frac{z+r{e}^{{i}\theta }}{1+r{e}^{-{i}\theta }z}=\frac{r{i}%
(1+2rz+z^{2})}{(1+rz)^{2}}\text{.}
\end{eqnarray*}

Next we compute 
\begin{equation*}
2\pi \Vert X\Vert ^{2}=\int_{{\mathbb{S}}^{1}}\left\vert \widetilde{X}\left(
z\right) \right\vert ^{2}\,{d}m(z)=\int_{{\mathbb{S}}^{1}}\left\vert {i}%
T_{r}(z)\right\vert ^{2}\,{d}m(z)=\int_{{\mathbb{S}}^{1}}1\,{d}m(z)=2\pi 
\text{.}
\end{equation*}%
We have also 
\begin{equation*}
2\pi \Vert Y\Vert ^{2}=\int_{{\mathbb{S}}^{1}}\left\vert \widetilde{Y}\left(
z\right) \right\vert ^{2}\,{d}m(z)=\int_{{\mathbb{S}}^{1}}\left\vert \frac{%
1-z^{2}}{(1+rz)^{2}}\right\vert ^{2}\,{d}m(z)\text{.}
\end{equation*}%
Setting $z={e}^{{i}s}$, we have%
\begin{equation*}
2\pi \Vert Y\Vert ^{2}=\int_{0}^{2\pi }\frac{1}{{i}{e}^{{i}s}}\left\vert 
\frac{1-{e}^{{i}s2}}{(1+r{e}^{{i}s})^{2}}\right\vert ^{2}{i}{e}^{{i}s}\,{d}%
s=\int_{{\mathbb{S}}^{1}}\frac{1}{{i}z}\left\vert \frac{1-z^{2}}{(1+rz)^{2}}%
\right\vert ^{2}{d}z\text{.}
\end{equation*}%
Now, the integrand is a complex analytic function inside the circle (observe
that $\bar{z}=1/z$ for $\left\vert z\right\vert =1$), except for a simple
pole at $z=0$ and a pole of order two at $z=-r$, with residues $\frac{{i}}{%
r^{2}}$ and $\frac{{i}(r^{2}+1)}{-r^{2}(1-r^{2})}$, respectively. One
obtains that $\Vert Y\Vert ^{2}=2/\left( 1-r^{2}\right) $. In the same way
one gets $\Vert Z\Vert ^{2}=2r^{2}/\left( 1-r^{2}\right) $.

We claim that the vectors $X,Y,Z$ are pairwise orthogonal. Let $h\left(
U,V\right) =U\bar{V}$ denote the Hermitian inner product on $\mathbb{C}$. We
compute 
\begin{equation*}
\int_{{\mathbb{S}}^{1}}h\left( \widetilde{X}\left( z\right) ,\widetilde{Y}%
\left( z\right) \right) \,{d}m(z)=\int_{{\mathbb{S}}^{1}}f(z)\,{d}z\text{,}
\end{equation*}%
where $f\left( z\right) =\frac{z^{2}-1}{z\left( 1+rz\right) \left(
z+r\right) }$ is a complex analytic function inside the circle, except for
simple poles at $z=0$ and $z=-r$, with residues $1/r$ and $-1/r$,
respectively. Then, 
\begin{equation*}
\left\langle X,Y\right\rangle =\Re \int_{{\mathbb{S}}^{1}}h\left( \widetilde{%
X}\left( z\right) ,\widetilde{Y}\left( z\right) \right) \,{d}m(z)=0\text{.}
\end{equation*}%
Analogously, we find that $\left\langle Y,Z\right\rangle =\left\langle
X,Z\right\rangle =0$. \hfill $\square $

\medskip

\begin{corollary}
The force free Moebius motions of the circle, or equivalently, the
geode\-sics of $\mathcal{M}$, are, via $F$, of the form $\gamma =\left(
\gamma _{1},\gamma _{2}\right) $, where $\gamma _{1}$ parametrizes the
circle with constant speed and $\gamma _{2}$ is a geodesic in the disc $%
\Delta $ whose trajectory coincides with the images of either $c_{1}(\rho
)=(\rho ,\theta _{0})$ or $c_{2}\left( \theta \right) =(\rho (\theta
),\theta )$, where $\rho $ satisfies the differential equation 
\begin{equation}
\left( \rho ^{\prime }\right) ^{2}=\frac{\mu +\rho ^{2}}{\left( 1-\rho
^{2}\right) \rho ^{2}}  \label{ecuacion}
\end{equation}%
for some constant $\mu >-1$.
\end{corollary}

\noindent Proof. Clearly, a geodesic of a Riemannian product projects to a
geodesic in each factor. Besides, as the coefficients of the first
fundamental form of $\Delta $ depend only on $\rho $, the corresponding
metric is Clairaut. Then, the trajectories of the geodesics of $\Delta $
are, in polar form, 
\begin{equation*}
c_{1}(\rho )=(\rho ,\theta _{0})\text{ \ \ \ \ or \ \ \ \ \ \ \ }c_{2}\left(
\theta \right) =(\rho (\theta ),\theta )
\end{equation*}%
for some constant $\theta _{0}$, where $\rho (\theta )$ satisfies Clairaut's
differential equation, for some $\lambda $: 
\begin{equation*}
\lambda E^{2}(\rho )=E(\rho )+(\rho ^{\prime })^{2}G(\rho )\text{.}
\end{equation*}%
Since in our case $E(\rho )=\left\Vert \frac{\partial }{\partial r}%
\right\Vert ^{2}=\frac{2}{1-\rho ^{2}}$ and $G(\rho )=\left\Vert \frac{%
\partial }{\partial \theta }\right\Vert ^{2}=\frac{2\rho ^{2}}{1-\rho ^{2}}$%
, the differential equation is equivalent to (\ref{ecuacion}) for some
constant $\mu >-1$. \hfill $\square $

\bigskip

\noindent Daniela Emmanuele

\noindent Departamento de Matem\'{a}tica, Escuela de Ciencias Exactas y
Naturales,

\noindent Facultad de Ciencias Exactas, Ingenier\'{\i}a y Agrimensura,

\noindent Universidad Nacional de Rosario, Argentina

\noindent emman@fceia.unr.edu.ar

\medskip

\noindent Marcos Salvai

\noindent CIEM - FaMAF, 

\noindent Conicet - Universidad Nacional de C\'ordoba, Argentina

\noindent salvai@famaf.unc.edu.ar

\end{document}